\input amstex
\documentstyle{amsppt}
\magnification=\magstep1 \pagewidth{6.5truein}
\pageheight{9.0truein} \loadbold

\TagsAsMath

\def\lmat{\left[\matrix}
\def\rmat{\endmatrix\right]}

\def\NN{{\Bbb N}}
\def\sr{\operatorname{sr}}
\def\mbull{M_\bullet(A)}
\def\End{\operatorname{End}}
\def\Hom{\operatorname{Hom}}

\def\almostiso{{\bf 1}}
\def\Bass{{\bf 2}}
\def\Blackstab{{\bf 3}}
\def\Blacksr{{\bf 4}}
\def\HV{{\bf 5}}
\def\Lamb{{\bf 6}}
\def\LamLect{{\bf 7}}
\def\Lam{{\bf 8}}
\def\Rie{{\bf 9}}
\def\Va{{\bf 10}}

\topmatter

\title Stable rank of corner rings  \endtitle

\author P. Ara and K. R. Goodearl
\endauthor

\rightheadtext{stable rank of corner rings}

\address {P. Ara: Departament de Matem\`atiques, Universitat
Aut\`onoma de Bar\-ce\-lo\-na, 08193 Bellaterra (Barcelona),
Spain}\endaddress \email para\@mat.uab.es\endemail

\address {K. R. Goodearl: Department of Mathematics, University of
California, Santa Barbara, CA 93106,
USA}\endaddress\email{goodearl\@math.ucsb.edu}\endemail

\abstract B.~Blackadar recently proved that any full corner $pAp$ in a
unital C*-algebra $A$ has K-theoretic stable rank greater than or
equal to the stable rank of $A$. (Here $p$ is a projection in $A$, and
fullness means that $ApA=A$.) This
result is extended to arbitrary (unital) rings $A$ in the present
paper: If $p$ is a full idempotent in $A$, then $\sr(pAp)\ge \sr(A)$.
The proofs rely partly on algebraic analogs of Blackadar's methods, and
partly on a new technique for reducing problems of higher stable rank
to a concept of stable rank one for skew (rectangular) corners $pAq$.
The main result yields estimates relating stable ranks of Morita
equivalent rings. In particular, if $B\cong \End_A(P)$ where $P_A$ is a
finitely generated projective generator, and $P$ can be generated by
$n$ elements, then $\sr(A)\le n{\cdot}\sr(B)-n+1$.
\endabstract

\thanks
 The first-named author was partially supported by the DGI and European
Regional Development Fund, jointly, through Project BFM2002-01390, and
by the Comissionat per Universitats i Recerca de la Generalitat de
Catalunya. The second-named author was partially supported by an NSF
grant. Part of this work was done during his research stay at the
Centre de Recerca Matem\`atica (Barcelona) in Spring 2003; he thanks the
CRM for its hospitality and support.
\endthanks



\endtopmatter

\document

\head Introduction \endhead

The theory of stable range of rings was developed by H. Bass
\cite{\Bass} and L. N. Vaserstein \cite{\Va}. As is now common, we
define the {\it stable rank\/} of a ring $A$, denoted $\sr(A)$, to be
the least positive integer $n$ such that $A$ satisfies Bass's $n$-th
stable range condition, or $\infty$ if no such $n$ exists. It is well
known that stable rank is not Morita invariant. In fact, Vaserstein
\cite{\Va} computed the stable rank of a matrix ring
$M_n(A)$, obtaining the following amazing formula
$$\sr(M_n(A))={\biggl\lceil \frac{\sr (A)-1}{n}\biggr\rceil}+1,$$
where $\lceil r\rceil$ denotes the least integer greater than or equal
to a real number $r$. If $B$ is a ring Morita equivalent to $A$, then
$B\cong pM_n(A)p$ for some full idempotent $p\in M_n(A)$. Thus, to
understand the behavior of stable rank under Morita equivalence, it
remains to see what happens to stable rank under the passage from a
ring to a full corner. Vaserstein's formula already contains some
information in this direction, namely that $\sr(A)\ge \sr(M_n(A))$ for
all $n\in \NN$. Since
$A$ is isomorphic to a corner ring in $M_n(A)$, corresponding to the
full idempotent $e_{11}$, this suggests the inequality $\sr(pAp)\ge
\sr(A)$ for any full corner $pAp$ of $A$. Such a formula was
conjectured by Blackadar \cite{\Blackstab, Remark A7} to hold for the
{\it topological stable rank\/} introduced by Rieffel in
\cite{\Rie}. It was subsequently proved by Herman and Vaserstein
\cite{\HV} that the Rieffel topological stable rank and the Bass stable
rank agree for any C*-algebra. Blackadar has recently verified
the corner conjecture in \cite{\Blacksr, Theorem 4.5}. His methods
are focussed on the topological stable rank, and rely on norm estimates
for differences of row vectors.

Previous work on stable rank of corners gave weaker inequalities of the
following form. If $p$ is a full projection in a unital C*-algebra $A$,
then, for some $n$, there exist $n$ pairs $(a_i,b_i)\in A^2$ such that
$\sum_{i=1}^n a_ipb_i =1$. Blackadar showed in \cite{\Blackstab, Lemma
A6} that in this situation, $\sr(A)\le \sr(pAp)+n-1$. That this result
extends to full idempotents in arbitrary rings was noted by the present
authors in
\cite{\almostiso, Remark 1.4}. In particular, this inequality suffices
to show that finiteness of the stable rank is Morita invariant.

Here we prove that the inequality $\sr(A)\le \sr (pAp)$
holds for any full corner $pAp$ in any unital ring $A$ (Theorem 7).
 The structure of the proof has been modelled after Blackadar's paper
\cite{\Blacksr}, but we have had to replace his topological
methods with purely algebraic ones. Of crucial importance is the
notion of stable rank one for skew, or rectangular, corners $pAq$,
where $p$ and $q$ are distinct idempotents of $A$. This allows us
to work only with stable rank one conditions, thus avoiding higher
rank conditions. By combining our main result with Vaserstein's
formula, we obtain estimates comparing the stable ranks of Morita
equivalent rings (Theorem 9).

We note that Lam and Dugas \cite{\Lam} have recently shown that
the reverse inequality $\sr(A)\ge \sr(eAe)$ holds for any
quasi-duo ring $A$ and any idempotent $e$ in $A$. By definition, a
{\it quasi-duo ring\/} is a ring in which every maximal one-sided
ideal is an ideal. It is clear that the only full idempotent in a
quasi-duo ring is $1$, so our result does not give any further
insight into Lam and Dugas's, nor vice versa.

\head Stable rank and skew corners \endhead

Throughout, let $A$ be a unital ring. We start by recalling the
definition of the (Bass) stable rank:

\definition{Definition} An $n$-row $(a_1,\dots ,a_n)\in A^n$ is said to
be {\it right unimodular\/} if $\sum _{i=1}^n a_iA=A$. An $(n+1)$-row
$(a_1,\dots ,a_n,b)\in A^{n+1}$ is {\it reducible\/} in case there
is an $n$-row $(c_1,\dots ,c_n)\in A^n$ such that the $n$-row
$(a_1+bc_1,\dots,a_n+bc_n)$ is right unimodular. The {\it stable
rank\/} of $A$, denoted $\sr (A)$, is the least positive integer $n$
such that every right unimodular $(n+1)$-row in $A^{n+1}$ is reducible,
or $\infty $ if no such $n$ exists.
\enddefinition

We next recall some
useful terminology. Two idempotents $p$ and $q$ in $A$ are said to
be {\it orthogonal\/}, written $p\perp q$, in case $pq=qp=0$. The
set of all idempotents of $A$ is partially ordered by declaring
$p\le q$ if and only if $p=pq=qp$. The idempotents $p$ and $q$ are
{\it equivalent\/}, written $p\sim q$, in case there are elements
$a\in pAq$ and $b\in qAp$ such that $p=ab$ and $q=ba$. (Note that
$p$ and $q$ are equivalent if and only if the right (respectively, left)
ideals generated by $p$ and $q$ are isomorphic as a right
(respectively, left) $R$-modules \cite{\Lamb, Proposition 21.20}.)
We write $p\lesssim q$ in case there is an idempotent $p'$ such
that $p'\le q$ and $p\sim p'$; this occurs if and only if there exist
elements $a\in pAq$ and $b\in qAp$ such that $ab=p$. For any idempotents
$p,q\in A$, we write $p\oplus q$ for the idempotent $\text{diag}(p,q)$
in
$M_2(A)$. Accordingly, the notation $n{\cdot}p$ is used for the
idempotent $\text{diag}(p,p, \dots ,p)$ in $M_n(A)$.

For all $m,n\in\NN$, identify $M_n(A)$ with the $n\times n$ upper
left corner subring of $M_{n+m}(A)$. In particular, $A= M_1(A)$ is
then identified with a subring of each $M_n(A)$. With this
identification, $1_A$ equals the matrix unit $e_{11}$, and $n{\cdot}1_A$
equals the identity matrix in $M_n(A)$. These identifications
allow us to work in as large a matrix ring as is convenient. When
the size of the matrices is not relevant, we write $\mbull$ to
stand for $M_n(A)$ with $n$ unspecified.

We say that an idempotent $p$ in $A$ is {\it full\/} in case $p$
generates $A$ as a two-sided ideal, that is, $ApA=A$. It is standard
(and an easy exercise) that
$p$ is a full idempotent if and only if $1_A\lesssim t{\cdot}p$ for some
$t\in \NN$. A {\it corner\/} of $A$ is any subring of the form $pAp$,
where $p$ is an idempotent, and we say that $pAp$ is a {\it full
corner\/} in case $p$ is a full idempotent. A {\it skew\/} ({\it
rectangular\/}) {\it corner\/} in $A$ is any subset of the form $pAq$,
for idempotents $p,q\in A$. Note that $pAq$ is a $(pAp,qAq)$-bimodule,
and that the ring multiplication in $A$ induces bimodule homomorphisms
$pAq\otimes_{qAq} qAp \rightarrow pAp$ and $qAp\otimes_{pAp}
pAq\rightarrow qAq$.

\definition{Definition} Let $p,q\in A$ be idempotents. We say that
the skew corner $pAq$ has (right) {\it stable rank $1$\/},
abbreviated $\sr(pAq) =1$, provided the following condition holds:
Whenever $a\in pAq$, $x\in qAp$, and $b\in pAp$ such that
$ax+b=p$, there exist $y\in pAq$ and $z\in qAp$ such that $(a+by)z
= p$. Note that in case $p=q$, there is no conflict between this
definition and the statement that the stable rank of the ring $pAp$ is
1. \enddefinition

The key to our methods is the following lemma, which reduces stable
rank calculations to questions of stable rank 1 for skew corners. Note
that $1_AM_n(A)$ is a skew corner, since it equals
$1_AM_n(A)(n{\cdot}1_A)$.

\proclaim{Lemma 1} Let $n\in\NN$. Then $\sr(A)\le n$ if and only
if $\sr\bigl( 1_AM_n(A) \bigr) =1$. \endproclaim

\demo{Proof} ($\Longrightarrow$): Let $\alpha\in 1_AM_n(A)$,
$\chi\in M_n(A)1_A$, and $\beta\in 1_AM_n(A)1_A$ such that
$\alpha\chi+ \beta= 1_A$. Then
 $$\xalignat3 \alpha &= \lmat a_1 &a_2 &\cdots &a_n\\ 0 &0 &\cdots &0\\
 \vdots &\vdots &&\vdots\\ 0 &0 &\cdots &0 \rmat,  &\chi &= \lmat x_1 &0
 &\cdots &0\\ x_2 &0
 &\cdots &0\\ \vdots &\vdots &&\vdots\\ x_n &0 &\cdots &0 \rmat,  &\beta
 &= \lmat b &0 &\cdots &0\\ 0 &0 &\cdots &0\\ \vdots &\vdots &&\vdots\\ 0
 &0 &\cdots &0 \rmat \tag\dagger \endxalignat$$
for some $a_i,x_i,b\in A$ such that $a_1x_1+ \cdots+ a_nx_n +b
=1$. Since $\sr(A)\le n$, there exist $y_1,\dots,y_n\in A$ such
that the row $(a_1+by_1, \dots, a_n+by_n)$ is right unimodular,
that is, $(a_1+by_1)z_1+ \cdots+ (a_n+by_n)z_n =1$ for some
$z_i\in A$. Setting
 $$\xalignat2 \zeta &= \lmat y_1 &y_2 &\cdots &y_n\\ 0 &0 &\cdots &0\\
 \vdots &\vdots &&\vdots\\ 0 &0 &\cdots &0 \rmat \in 1_AM_n(A),  &\xi &=
 \lmat z_1 &0 &\cdots &0\\ z_2 &0 &\cdots &0\\ \vdots &\vdots &&\vdots\\
 z_n &0 &\cdots &0 \rmat \in
 M_n(A)1_A, \tag\ddagger \endxalignat$$
we have $(\alpha+ \beta\zeta) \xi= 1_A$.

($\Longleftarrow$): If $(a_1,\dots,a_n,b)\in A^{n+1}$ is a right
unimodular row, there exist $x_1,\dots,x_n,x\in A$ such that
$a_1x_1+ \cdots+ a_nx_n+ bx= 1$. After replacing $b$ by $bx$, we
may assume that $x=1$. Define matrices $\alpha\in 1_AM_n(A)$,
$\chi\in M_n(A)1_A$, and $\beta\in 1_AM_n(A)1_A$ as in $(\dagger)$
above, and observe that $\alpha\chi+ \beta= 1_A$. Since
$\sr \bigl( 1_AM_n(A) \bigr) =1$, there exist matrices
$\zeta$ and $\xi$ as in $(\ddagger)$ such that $(\alpha+ \beta\zeta)
\xi= 1_A$. It follows that the row
$(a_1+by_1, \dots, a_n+by_n)$ is right unimodular. \qed\enddemo

We next show how stable rank 1 can be transferred from certain skew
corners to others.

\proclaim{Lemma 2} Let $p,q\in A$ be idempotents.

{\rm (a)} If $\sr(pAq) =1$, then $p\lesssim q$.

{\rm (b)} If $p',q'\in A$ are idempotents such that $p'\sim p$ and
$q'\sim q$, then $\sr(pAq) =1$ if and only if $\sr(p'Aq') =1$.
\endproclaim

\demo{Proof} (a) Consider the equation $0{\cdot}0+p =p$, where we
view the first $0\in pAq$, the second $0\in qAp$, and $p\in pAp$.
The hypothesis $\sr(pAq) =1$ then gives us $y\in pAq$ and $z\in
qAp$ such that $(0+py)z =p$. Hence, $pyqzp =p$, and it follows
that $p\lesssim q$.

(b) There are elements $u\in pAp'$ and $u'\in p'Ap$ such that
$uu'=p$ and $u'u=p'$, and elements $v\in qAq'$ and $v'\in q'Aq$
such that $vv'=q$ and $v'v=q'$.

Assume that $\sr(pAq)=1$, and consider elements $a\in p'Aq'$,
$x\in q'Ap'$, and $b\in p'Ap'$ such that $ax+b =p'$. Then we have
$uav'\in pAq$, $vxu'\in qAp$, and $ubu'\in pAp$ satisfying
$(uav')(vxu')+ (ubu')= p$. Since $\sr(pAq)=1$, there exist $y\in
pAq$ and $z\in qAp$ such that $(uav'+ ubu'y)z =p$. Then $u'yv\in
p'Aq'$ and $v'zu\in q'Ap'$ are elements satisfying the equation
$(a+ bu'yv) v'zu= p'$, which proves that $\sr(p'Aq') =1$. The
converse follows by symmetry. \qed\enddemo

\proclaim{Lemma 3} Let $p,q,s\in A$ be idempotents such that
$s\perp q$. If $\sr(pAq)= 1$, then also $\sr \bigl( pA(q+s) \bigr) =1$.
\endproclaim

\demo{Proof} Let $a\in pA(q+s)$, $x\in (q+s)Ap$, and $b\in pAp$
such that $ax+b =p$. Rewrite this equation as $(aq)(qx)+ (b+asx)
=p$, where $aq\in pAq$, $qx\in qAp$, and $b+asx\in pAp$. Since
$\sr(pAq) =1$, there exist $y\in pAq$ and $z\in qAp$ such that
 $$\bigl[ (aq)+ (b+asx)y \bigr] z =p.$$
Note that $sz=0$ because $s\perp q$, which allows us to rewrite the
equation above as
 $$\bigl[ a(q+s+sxy) +by \bigr] z =p.$$
Since $ys=0$, the element $sxy$ is nilpotent, and so the element $u
:=q+s+sxy$ is a unit in the ring $(q+s)A(q+s)$. Now
 $$(a+byu^{-1})(uz)= (au+by)z= p$$
with $yu^{-1}\in pA(q+s)$ and $uz\in (q+s)Ap$, which completes the
proof. \qed\enddemo

\proclaim{Lemma 4} Let $p,q,r\in A$ be idempotents such that
$p,q\perp r$. If $\sr \bigl( (p+r)A(q+r) \bigr) =1$, then $\sr(pAq) =1$.
\endproclaim

\demo{Proof} Let $a\in pAq$, $x\in qAp$, and $b\in pAp$ such that
$ax+b =p$. Then the elements $a+r\in (p+r)A(q+r)$ and $x+r\in
(q+r)A(p+r)$ satisfy $(a+r)(x+r) +b =p+r$. Since $\sr \bigl(
(p+r)A(q+r) \bigr) =1$, there exist $y\in (p+r)A(q+r)$ and $z\in
(q+r)A(p+r)$ such that
$(a+r+by)z =p+r$. Observe that
 $$\xalignat2 rz &= r(a+r+by)z =r  &(a+by)z &= p(a+r+by)z =p.
 \endxalignat$$
Then $z= (q+r)z= qz+r$, and so $zp= qzp\in qAp$. Since
 $$(a+bpyq)(zp)= (a+by)zp= p,$$
we have shown that $\sr(pAq) =1$. \qed\enddemo

\head The main results\endhead

The final ingredient needed to prove our main theorem is a (partial)
converse to Lemma 4, which holds when $r\in ApA$. The following
observations will be helpful.

\definition{Observation 5} If we are trying to establish
$\sr(pAq)=1$ for some idempotents $p,q\in A$, then we are given
$ax+b =p$ for some $a\in pAq$, $x\in qAp$, and $b\in pAp$, and we
seek $y\in pAq$ and $z\in qAp$ such that $(a+by)z =p$. Several
reduction steps are possible, in which we may do any of the
following:
 \roster
 \item Replace $a$ by $a+bc$ for any $c\in pAq$;
 \item Replace $a$ by $au$ for any unit $u$ of $qAq$;
 \item Replace $a$ by $va$ for any unit $v$ of $pAp$;
 \item Replace $A$ by $M_n(A)$ for any $n\in\NN$.
 \endroster
In cases (1)--(3), the replacement of $a$ by another element of $pAq$
must be accompanied by corresponding replacements for $x$ and $b$.

To see why (1) is allowed, for instance, observe that
 $$(a+bc)x +b(p-cx) =p$$
with $a+bc\in pAq$ and $b(p-cx)\in pAp$; if there exist $y'\in
pAq$ and $z'\in qAp$ such that
 $$\bigl[ (a+bc)+ b(p-cx)y' \bigr] z' =p,$$
then $\bigl[ a+ b(c+y'-cxy') \bigr] z' =p,$ with $c+y'-cxy'\in
pAq$. For (2), we have $(au)(u^{-1}x) +b =p$, and if $(au+by') z'
=p$, then $(a+by'u^{-1}) (uz') =p$. In the case of (3),
we have $(va) (xv^{-1}) +(vbv^{-1}) =p$, and if $(va+ vbv^{-1}y')
z' =p$, then $(a+bv^{-1}y') (z'v) =p$.

Finally, we address (4). Because of our identification of $A$ with the
corner $e_{11}M_n(A)e_{11}$, we have $pAq= pM_n(A)q$, and similarly for
$qAp$ and $pAp$. Thus, if there exist $y\in pM_n(A)q$ and $z\in
qM_n(A)p$ satisfying $(a+by)z =p$, then $y\in pAq$ and $z\in qAp$, and
the equation holds in $A$.
\enddefinition

 \proclaim{Proposition 6} Let $p,q,r\in A$ be
idempotents such that $p,q\perp r$. If $\sr(pAq) =1$ and $r\lesssim
n{\cdot}p$ for some $n\in\NN$, then $\sr \bigl( (p+r)A(q+r) \bigr) =1$.
\endproclaim

\demo{Proof} By Lemma 2, $p\sim p'$ for some idempotent $p'\le q$,
and we may replace $p$ by $p'$. Thus, we may assume that $p\le q$.

We claim that it suffices to establish the proposition under the
additional hypothesis $r\sim p$. Given that case, it follows by
induction that $\sr \bigl( (2^m{\cdot}p)\mbull (q\oplus
(2^m-1){\cdot}p) \bigr)=1$ for all
$m\in\NN$. In the general case, we choose $m$ large enough that
$2^m-1\ge n$, so that $(2^m-1){\cdot}p= r'\oplus r''$ for some
orthogonal idempotents $r'$ and $r''$ with $r'\sim r$. In a lower right
corner of a suitably large matrix ring, we can find an idempotent $s$
such that $s\perp 1_A$ and $s\sim r''$. Thus, $p,q,r\perp s$ and
$(2^m-1){\cdot}p\sim r+s$. Consequently, $\sr \bigl( (p+r+s)\mbull
(q+r+s) \bigr) =1$, and the desired result follows from Lemma 4.
Therefore the claim holds, and we may assume that $p\sim r$. In
particular,
$\sr(rAq) =1$.

To prove that $\sr \bigl( (p+r)A(q+r) \bigr) =1$, we may work within
the ring
$(q+r)A(q+r)$. Hence, we may assume, for convenience, that $q+r
=1$. Now elements $a\in A$ can be viewed as formal matrices of the
form
 $$\lmat qaq &qar\\ raq &rar \rmat.$$
We mimic the proof that stable rank 1 passes from a ring to its
$2\times2$ matrix ring, which would be exactly our present
situation in case $p=q$. In order to allow for the possibility
that $p<q$, we must be careful to modify our matrices starting at
the lower right corner (where the entries come from $rAr$), rather
than starting at the upper left corner, where we can control $pAq$
but not $qAq$.

Let $a\in (p+r)A$, $x\in A(p+r)$, and $b\in (p+r)A(p+r)$ such that
$ax+b =p+r$. Note that $qaq= paq$ and $qar= par$. Now
 $$(raq)(qxr) +(rarxr+ rbr)= r(ax+b)r= r.$$
Since $\sr(rAq) =1$, there exist $y_1\in rAq$ and $z_1\in qAr$
such that
 $$\bigl[ (raq)+ (rarxr+ rbr)y_1 \bigr] z_1= r.$$
The factor $q$ in $raq$ can be dropped from the last equation
because $qz_1= z_1$. Hence,
 $$r \bigl[ a(1+rxy_1)+ by_1 \bigr] z_1 =r.$$
Since $y_1r =0$, the element $1+rxy_1$ is a unit in $A$. In view
of Observation 5, we may replace $a$ by $a(1+rxy_1)+ by_1$. Thus,
we may now assume that there exists $z_1\in qAr$ such that $raz_1
=r$.

Next, we replace $a$ by the element
 $$\lmat paq &par\\ raq &rar \rmat \lmat q &z_1(r-rar)\\ 0 &r
 \rmat = \lmat paq &*\\ raq &r \rmat,$$
which is allowed since $\lmat q &z_1(r-rar)\\ 0 &r \rmat$ is a
unit in $A$. At this stage, we have $rar=r$. After replacing $a$
by
 $$\lmat p &-par\\ 0 &r \rmat \lmat paq &par\\ raq &r \rmat =
 \lmat * &0\\ raq &r \rmat$$
(note that $\lmat p &-par\\ 0 &r \rmat$ is a unit in
$(p+r)A(p+r)$), we may assume in addition that $par =0$.

Now return to the equation $ax+b =p+r$, and observe that
 $$(paq)(qxp) +(pbp) = p(ax+b)p =p,$$
because $pa= paq$. Since $\sr(pAq) =1$, there exist $y_2\in pAq$
and $z_2\in qAp$ such that $(paq+ pbpy_2)z_2 =p$. Consequently,
$p(a+by_2)z_2 =p$. Note that $(a+by_2)r= ar$, and so neither of
the conditions $rar=r$ and $par=0$ is lost on replacing $a$ by
$a+by_2$. Thus, we may assume that there exists $z_2\in qAp$ with
$paz_2 =p$, whence
 $$\lmat paq &0\\ raq &r \rmat \lmat z_2 &0\\ -raz_2 &r \rmat =
 \lmat p &0\\ 0 &r \rmat.$$
In other words, we have found an element $z=\lmat z_2 &0\\ -raz_2
&r \rmat$ in $A(p+r)$ such that $az= p+r$, and therefore we have
proved that $(p+r)A$ has stable rank 1. \qed\enddemo

\proclaim{Theorem 7} If $p$ is a full idempotent in $A$, then
$\sr(A) \le \sr(pAp)$. \endproclaim

\demo{Proof} Assume that $\sr(pAp) =n< \infty$. By Lemma 1,
$pM_n(pAp)= pM_n(A)(n{\cdot}p)$ has stable rank 1.
Since $p$ is full, $1_A\lesssim t{\cdot}p$ for some $t\in\NN$.
Working in a suitably large matrix ring $R=\mbull$, we have
$\sr(pR(n{\cdot}p)) =1$ and we have room for an idempotent $r$
which is equivalent to $(t-1){\cdot}p$ and orthogonal to both $p$
and $n{\cdot}p$. Proposition 6 now implies that
$\sr \bigl( (p+r)R(n{\cdot}p+r) \bigr) =1$, that is,
$\sr \bigl( (t{\cdot}p)R((n+t-1){\cdot}p) \bigr) =1$. We also have
$t{\cdot}p= e+f$ and $(n+t-1){\cdot}p= e+f+g$ for some orthogonal
idempotents $e$, $f$, $g$ with
$e\sim 1_A$ and $g\sim (n-1){\cdot}p$. Since $\sr \bigl( (e+f)R(e+
f+g) \bigr) =1$, we can use Lemma 4 to see that
$\sr \bigl( eR(e+g) \bigr) =1$. Finally, using Lemma 3
to increase
$e+g$ by an orthogonal idempotent equivalent to
$(n-1){\cdot}(1_A-p)$, we conclude that
$\sr \bigl( eR(e\oplus (n-1){\cdot}1_A) \bigr) =1$. Since $e\sim 1_A$,
we thus have $\sr \bigl( 1_AR(n{\cdot}1_A) \bigr)=1$, and so
$\sr \bigl( 1_AM_{n}(A) \bigr) =1$. Therefore $\sr(A) \le n$, by
Lemma 1. \qed\enddemo

An upper bound for $\sr(pAp)$ in terms of $\sr(A)$ can be obtained from
Theorem 7 and Vaserstein's formula, as follows.

\proclaim{Theorem 8} If $p$ is a full idempotent in $A$ and
$\sum_{i=1}^n a_ipb_i =1$ for some $a_i,b_i\in A$, then $\sr(pAp)
\le n{\cdot}\sr(A) -n+1$. \endproclaim

\demo{Proof} We may clearly assume that each $a_i\in Ap$ and each
$b_i\in pA$. Set
 $$\align \alpha &= \lmat a_1 &a_2 &\cdots &a_n\\ 0 &0 &\cdots &0\\
 \vdots &\vdots &&\vdots\\ 0 &0 &\cdots &0 \rmat \in
1_AM_n(A)(n{\cdot}p) \\
\beta &=
 \lmat b_1 &0 &\cdots &0\\ b_2 &0 &\cdots &0\\ \vdots &\vdots &&\vdots\\
 b_n &0 &\cdots &0 \rmat \in
(n{\cdot}p)M_n(A)1_A, \endalign$$ and observe that $\alpha\beta= 1_A$.
Then the matrix $q:= \beta\alpha$ is an idempotent in the ring
$(n{\cdot}p)M_n(A)(n{\cdot}p)$, which we identify with $M_n(pAp)$. Since
$(n{\cdot}p)\alpha q\beta(n{\cdot}p)= p$, we see that $q$ is full in
$M_n(pAp)$. Moreover, $1_A\sim q$ in $M_n(A)$, and so
$$A\cong 1_AM_n(A)1_A\cong qM_n(A)q= qM_n(pAp)q.$$
Hence, $\sr(M_n(pAp)) \le
\sr(A)$ by Theorem 7. According to Vaserstein's formula \cite{\Va,
Theorem 3},
 $$\sr(M_n(pAp))\ge \dfrac{\sr(pAp)-1}{n}+1,$$
and the theorem follows. \qed\enddemo

To conclude, we derive the following
estimates for the stable ranks of Morita equivalent rings.

\proclaim{Theorem 9} Let $A$ and $B$ be Morita equivalent rings; then
$B\cong \End_A(P)$ for some finitely generated projective
generator $P_A$. If $P$ can be generated by $n$ elements as a right
$A$-mod\-ule, then
$$\sr(A) \le n{\cdot}\sr(B)-n+1.$$
If there are $t$ homomorphisms $f_i\in \Hom_A(P,A)$ such
that $\sum_{i=1}^t f_i(P)=A$, then
$$\sr(B) \le t{\cdot}\sr(A)-t+1.$$
\endproclaim

\demo{Proof} There exists a split epimorphism $A^n
\twoheadrightarrow P$, so that $P\cong pA^n$ for some idempotent
$p\in M_n(A)$, and $B\cong pM_n(A)p$. Since $P$ is a generator, $p$ is
full. Thus, by Vaserstein's formula and Theorem 7,
$$\dfrac{\sr(A)-1}{n}+1\le \sr(M_n(A))\le \sr(B).$$

We now identify $B$ with $\End_A(P)$, and view $P$ as a
left $B$-module. Then $_BP$ is a finitely generated projective
generator, and $\End_B(P)\cong A$ (e.g., \cite{\LamLect, Propositions
18.17, 18.22}). There exist elements $x_i\in P$ such that $\sum_{i=1}^t
f_i(x_i) =1$ in
$A$. Given any element $x\in P$, there are endomorphisms $xf_i\in B$
(sending any
$y\mapsto xf_i(y)$) such that $\sum_{i=1}^t (xf_i)x_i= x\sum_{i=1}^t
f_i(x_i) =x$. This shows that $P$ is generated as a left $B$-module by
$x_1,\dots,x_t$. Therefore the inequality $\sr(B) \le
t{\cdot}\sr(A)-t+1$ follows from the first part of the theorem, on
replacing $A$ and $B$ by $B^{\text{op}}$ and $A^{\text{op}}$,
respectively.
\qed\enddemo

\enddocument